\documentclass[12pt]{article}

\usepackage{amsfonts,amsmath}
\usepackage{amssymb,latexsym}
\usepackage[dvips]{epsfig}
\usepackage{amscd}  

\title{\bf An invariant of tangle cobordisms}
\author{Mikhail Khovanov\footnote{Electronic address: 
  mikhail@math.ucdavis.edu}} 
\date{July 30, 2002}

\begin{document}
\maketitle
\baselineskip 14pt 

\def\R{\mathbb R}
\def\N{\mathbb N}
\def\Z{\mathbb Z}
\def\Q{\mathbb Q}
 
\def\l{\lbrace}
\def\r{\rbrace}
\def\lra{\longrightarrow}

\newcommand{\mc}{\mathcal} 
\newcommand{\mf}{\mathfrak}

\newcommand{\cA}{{\mathcal{A}}}   
\newcommand{\cH}{{\mathcal{H}}}   
\newcommand{\cC}{{\mathcal{C}}}   
\newcommand{\cF}{{\mathcal{F}}}   
\newcommand{\cK}{\mathcal{K}}     

\newcommand{\drawing}[1]{\begin{center} \epsfig{file=#1} \end{center}}

\newtheorem{prop}{Proposition}
\newtheorem{theorem}{Theorem}
\newtheorem{lemma}{Lemma}
\newtheorem{corollary}{Corollary}

\newcommand{\oplusop}[1]{{\mathop{\oplus}\limits_{#1}}}

\def\Hom{\mathrm{Hom}}
\newcommand\hsm{\hspace{0.05in}}
\newcommand\vsp{\vspace{0.1in}}
\newcommand{\Id}{\mathrm{Id}}
\newcommand{\define}{\stackrel{\mbox{\scriptsize{def}}}{=}}
\newcommand{\Vertical}{\mathrm{Vert}}

\section{Introduction} 

In \cite{me:tangles} to a plane diagram $D$ of an oriented tangle $T$
with $2n$ bottom and $2m$ top endpoints we associated a complex $\cF(D)$ 
of $(H^m,H^n)$-bimodules, for certain rings $H^n.$ 
We proved that the isomorphism class of this 
complex in the homotopy category is an invariant of $T.$ 
In this paper we give a short argument that 
 our construction yields an invariant of tangle 
cobordisms. To a diagram of an  oriented cobordism  between diagrams 
$D_1$ and $D_2$ of tangles 
$T_1$ and $T_2$ we assign a homomorphism of complexes 
$\cF(D_1)\to \cF(D_2)$ and then check that (in the homotopy 
category) this homomorphism depends on the choice of a diagram of the 
cobordism only up to the overall minus sign. The result follows 
from the basic properties of rings $H^n$ and $H^n$-bimodules
assigned to tangle diagrams. 

For link cobordisms this result was recently obtained by Magnus Jacobsson 
\cite{Jacobsson}. 

In a previous paper \cite{me:jones} we conjectured that such an invariant 
exists over the ring $\Z[c].$ This conjecture, which should be understood 
"rel boundary" (as emphasized by Jacobsson \cite{Jacobsson}), remains 
open. Jacobsson established the $c=0$ specialization, and it also follows 
from this work.

\section{2-tangles} 

The analogue of Reidemeister moves for surfaces embedded in $\R^4$ 
was found by Roseman \cite{Roseman} and investigated in depth by 
Carter and Saito  \cite{CS}, \cite{CSbook}. 

The framework for studying 2-tangles was developed  
by Fischer \cite{Fischer}, Kharlamov and Turaev \cite{KharTuraev},
Carter, Rieger, 
and Saito \cite{CRS}, and Baez and Langford \cite{BaezLang}. 
We will use a combinatorial 
realization of the 2-tangle 2-category described in 
\cite{CRS}, \cite[Section 2.5]{CSbook}, and \cite[Section 3]{BaezLang}. 
We assume familiarity with \cite{BaezLang}.  
Baez and Langford \cite{BaezLang} work with unoriented 2-tangles, but 
combinatorial description can be easily modified to the oriented case.  
We briefly review this description, referring the reader to 
\cite{BaezLang} for details. 

We consider oriented unframed tangles with even number of bottom  
endpoints and oriented cobordisms between these tangles. 

The objects of the 2-category $\mc{C}$ are even length sequences $s$ of 
pluses and minuses (indicating orientations of tangles near endpoints).  
Let $|s|$ denote half the length of $s.$ 

\begin{figure} [htb] \drawing{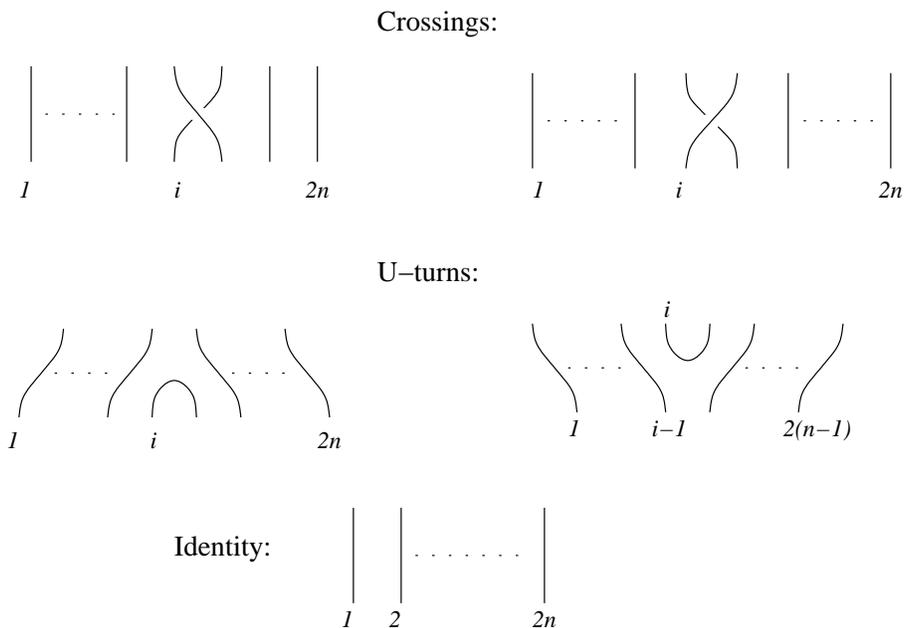}\caption{Generating 1-morphisms 
 of $\mc{C}$} 
 \label{gen1mor} 
 \end{figure}

1-morphisms of $\mc{C}$ represent planar diagrams of generic tangles. 
The generating 1-morphisms are positive and negative crossings, 
U-turns, and the identity 1-morphisms. They are depicted in 
figure~\ref{gen1mor}. Different orientations of arcs in a diagram lead to 
different 1-morphisms. We denote  U-turns by $\cap_{i,n}$ and $\cup_{i,n-1}$ 
and identity morphisms by $\Vertical_{n}$ (in \cite{me:tangles} we 
used $\Vertical_{2n}$ instead). 

 1-morphisms are products of generating 1-morphisms. Orientations of 
 arcs should be compatible when the diagrams are concatenated. 

\begin{figure} [htb] \drawing{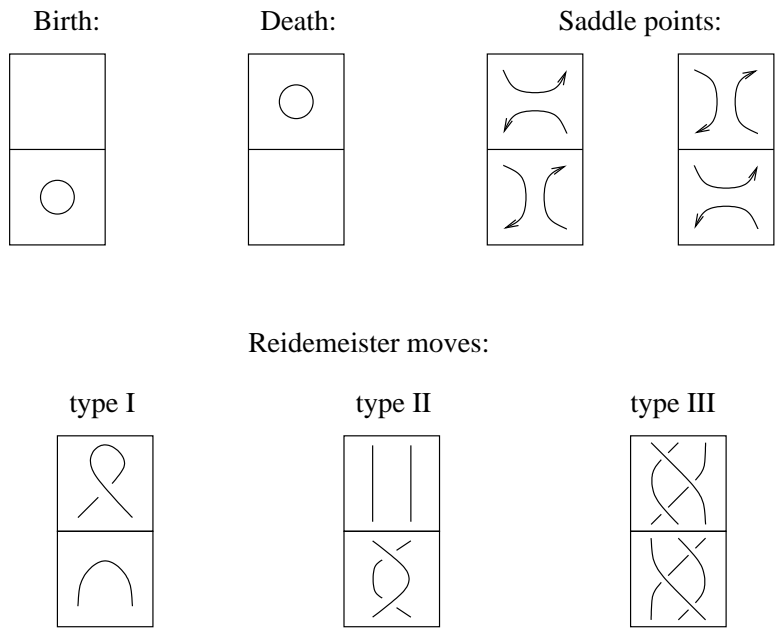}\caption{Generating 2-morphisms} 
 \label{2mor1} 
 \end{figure}

\begin{figure} [htb] \drawing{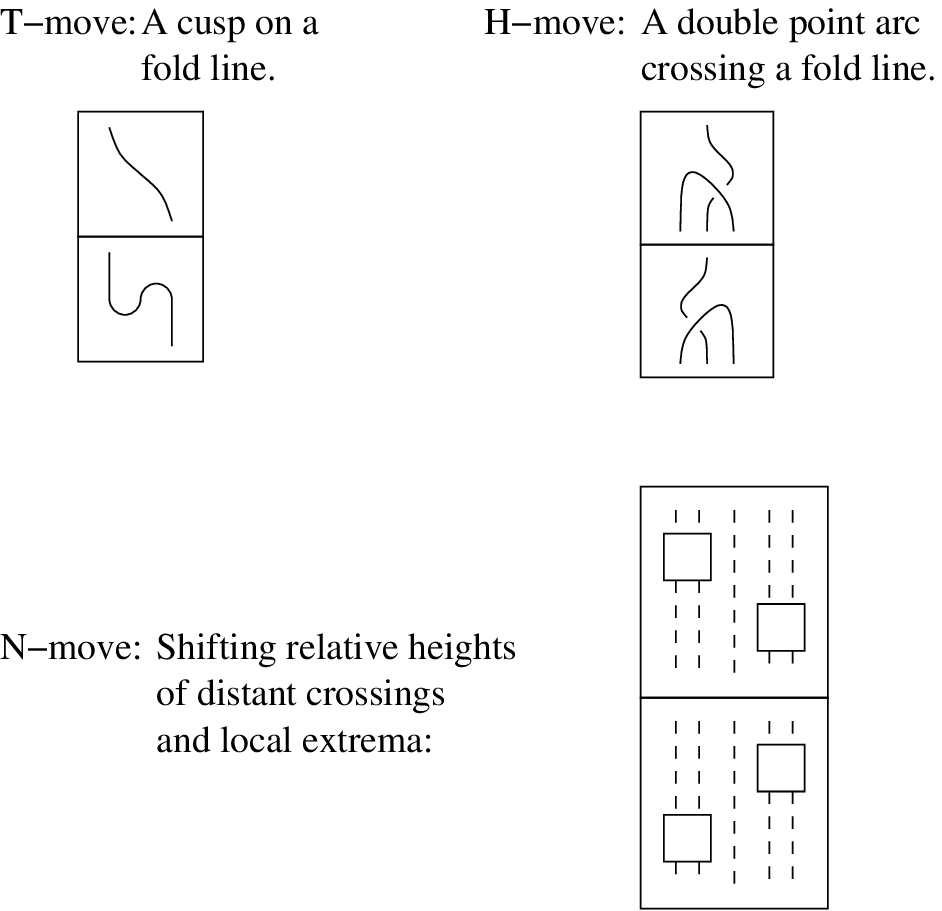}\caption{Generating 2-morphisms} 
 \label{2mor2} 
 \end{figure}

2-morphisms are combinatorial diagrams of tangle cobordisms, and 
depicted by "movies" of Roseman and Carter-Saito. 
The generating 2-morphisms are birth and death of a circle, 
saddle point (with compatible orientations), Reidemeister moves, 
a double point arc crossing a fold line, a cusp on a fold line, 
shifting relative heights of distant crossings and local extrema, 
and identity 2-morphisms. Generating 2-morphisms (except for 
identity morphisms) are depicted in figures~\ref{2mor1}, \ref{2mor2}. 

Each generating 2-morphism has several versions, obtained by 

(a) reading the film from bottom to top rather than from top to bottom, 

(b) changing between positive and negative crossings (the third Reidemeister 
 move has many such versions), 

(c) reflecting each frame about the x-axis,

(d) reflecting each frame about the y-axis, 

(e) orienting strings in various ways.   

Of course, for some moves some of these operations produce identical 
moves (and, for instance, operation (a) on a birth move produces 
a death move).

The figure~\ref{2mor2} two-morphisms will be called T-move, H-move, and
N-move, since these moves were labelled by T, H, and N in \cite{BaezLang} 
(with subscripts which we omit).   

\begin{figure}[htb] \drawing{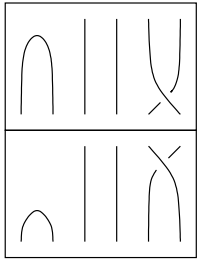}\caption{An example of height 
 shifting} 
 \label{exampleh} 
 \end{figure}

The height shifting morphism (N-move) has many versions, as we are free to put 
a U-turn or a crossing inside each small square of the top frame, 
add any number of strings separating the two small squares, and 
possibly invert the order of the film. 
An example is given in figure~\ref{exampleh}. 

\begin{figure} [htb] \drawing{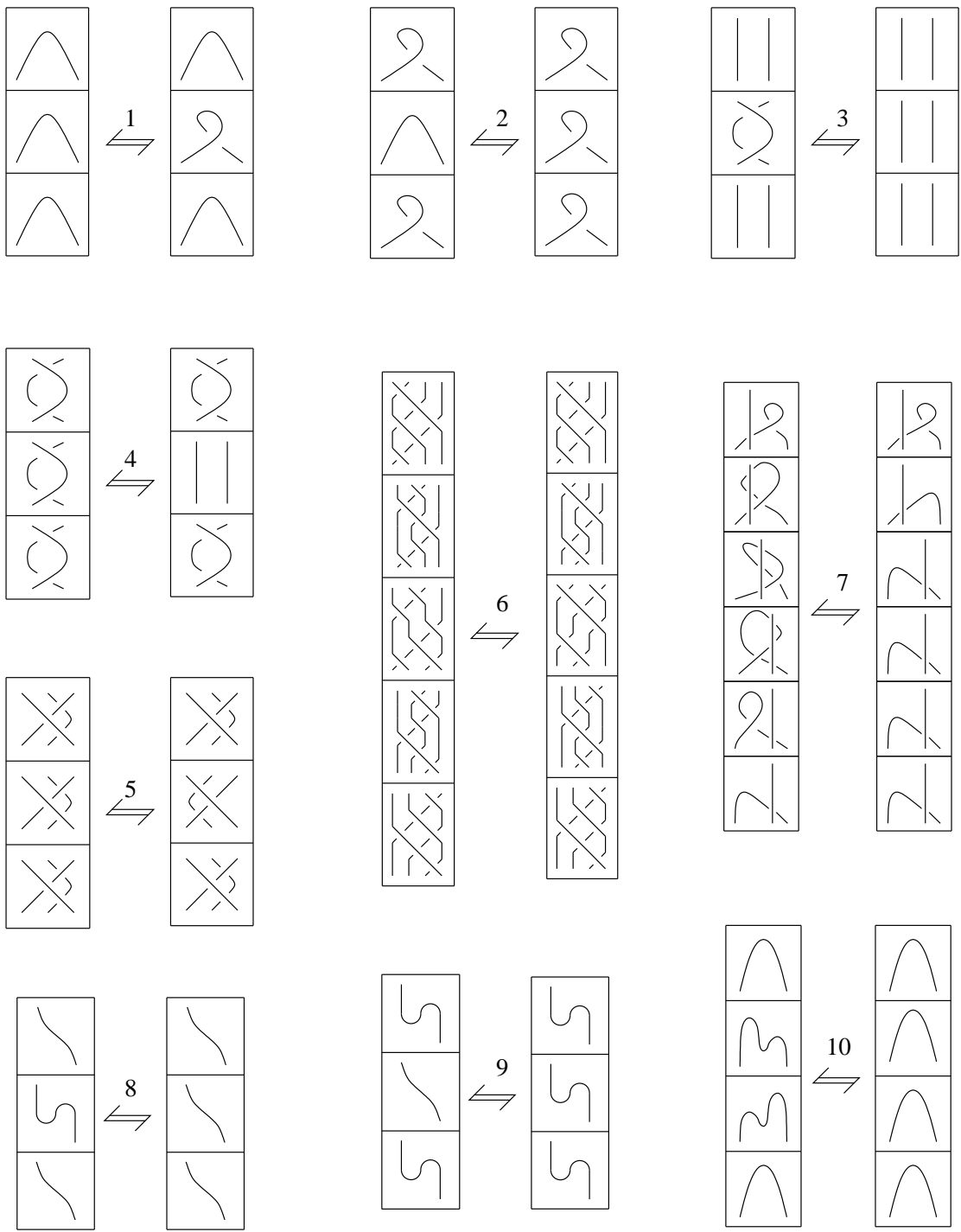}\caption{Movie moves 1-10} 
 \label{mm1} 
 \end{figure}

\begin{figure} [htb] \drawing{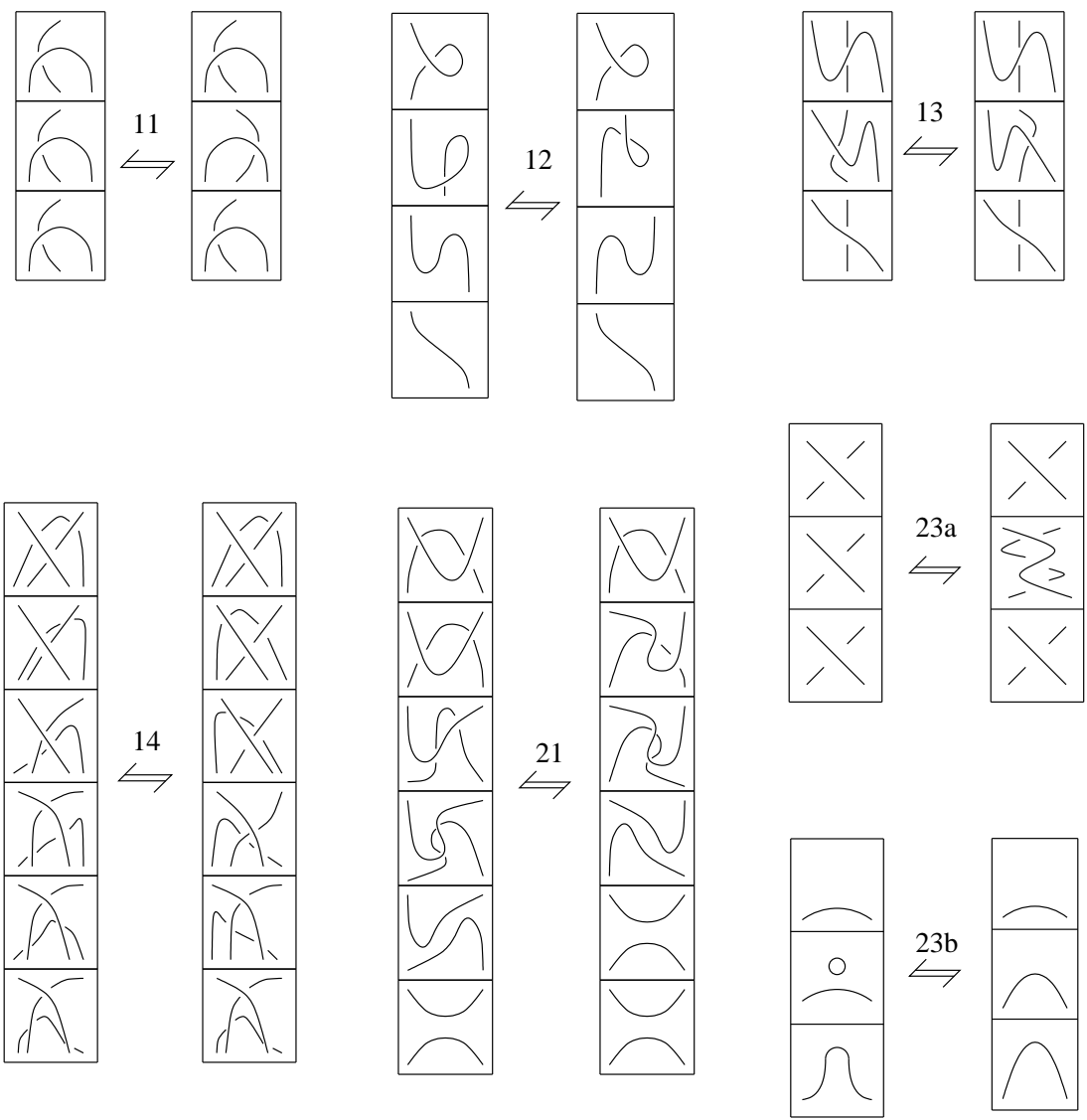}\caption{Movie moves 11-14, 21, 23} 
 \label{mm2} 
 \end{figure}

\begin{figure} [htb] \drawing{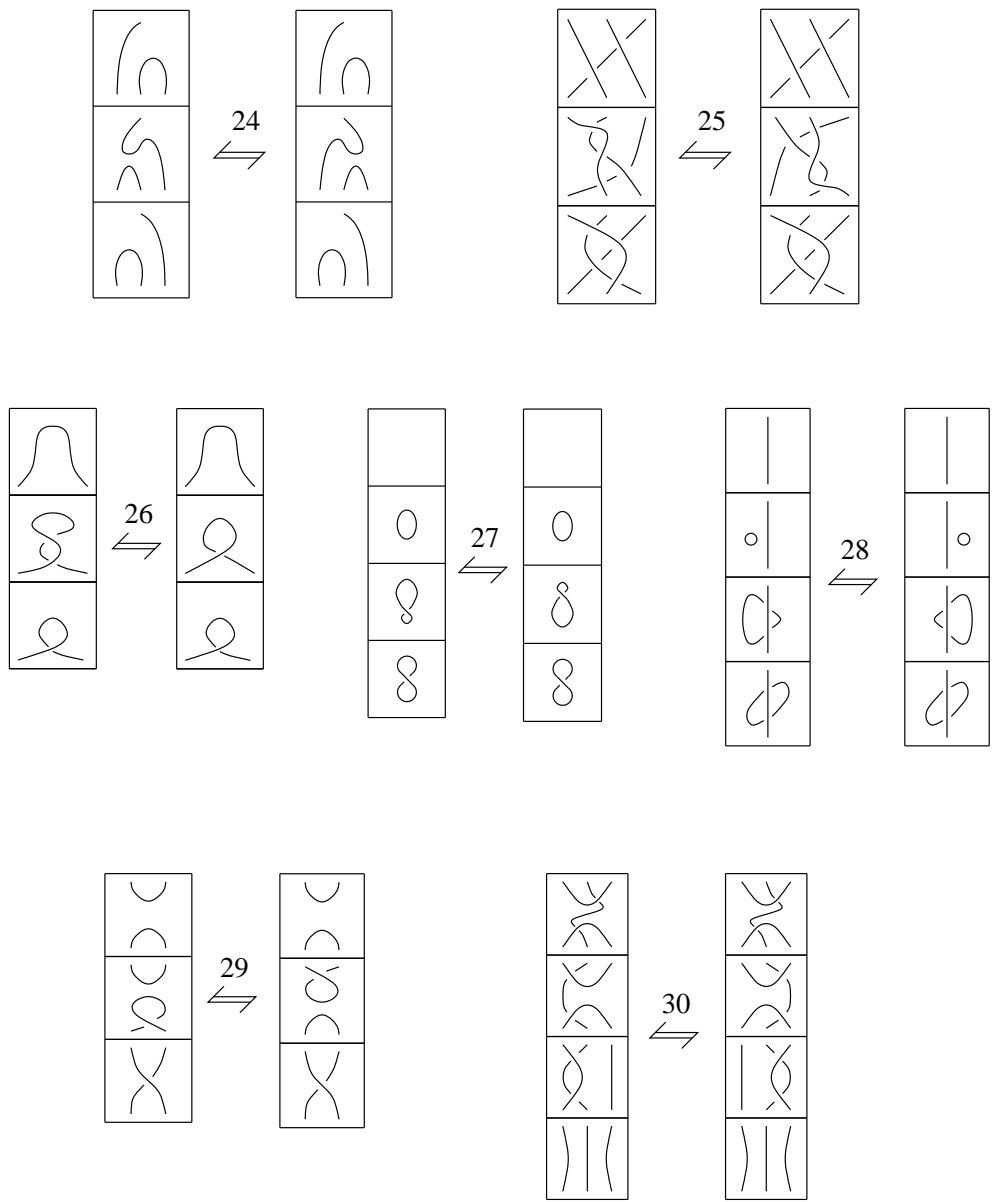}\caption{Movie moves 24-30} 
 \label{mm3} 
 \end{figure}

\begin{figure} [htb] \drawing{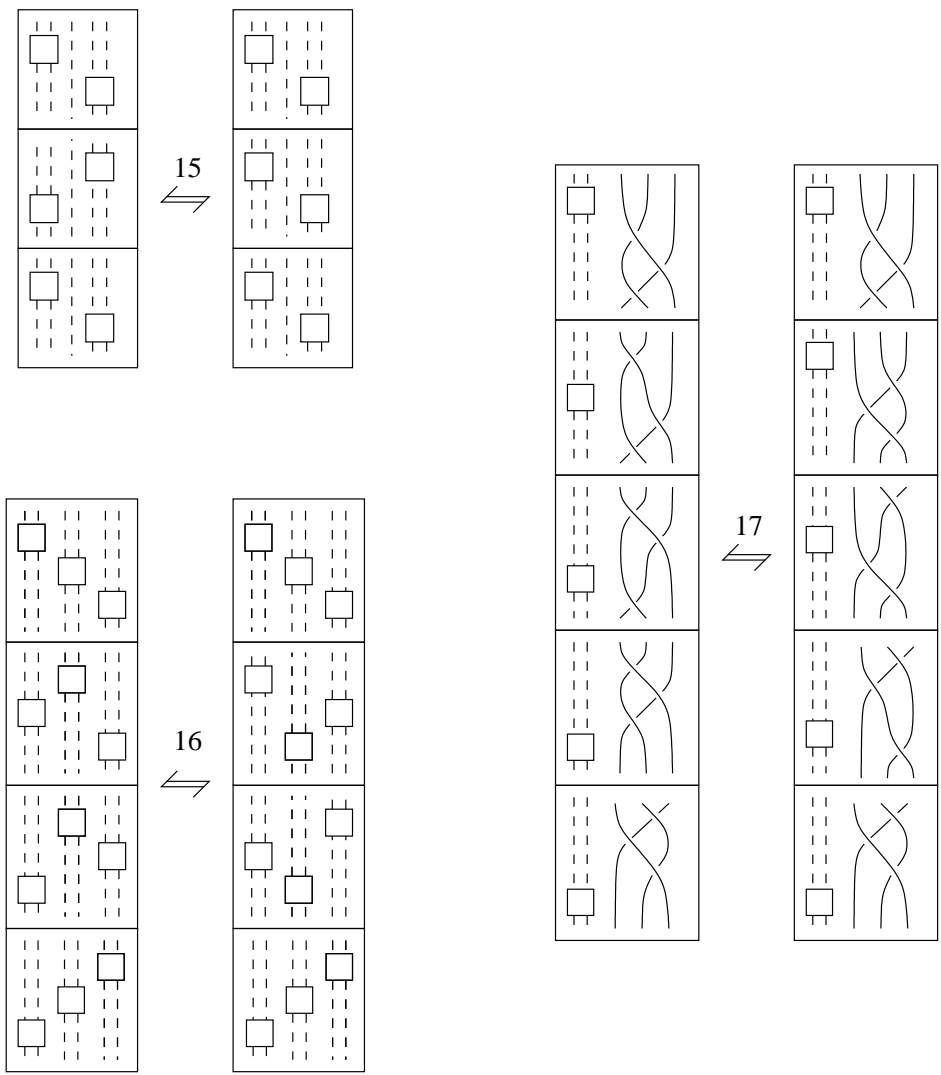}\caption{Movie moves 15-17} 
 \label{mm4} 
 \end{figure}

\begin{figure} [htb] \drawing{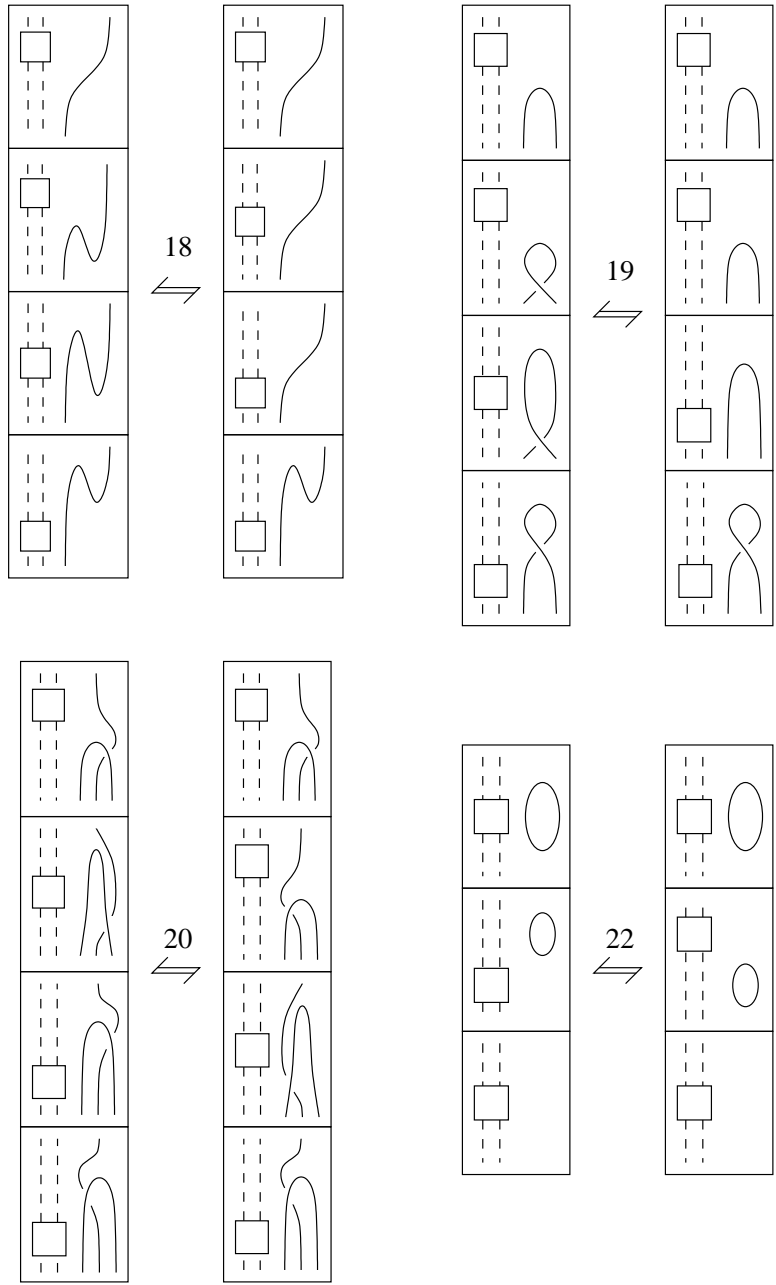}\caption{Movie moves 18-20, 22} 
 \label{mm5} 
 \end{figure}

A complete set of defining relations on 2-morphisms is given by 
the 31 movie moves (see \cite{CRS}, \cite[Section 2.5]{CSbook}, or 
  \cite{BaezLang}). The first 30 of these moves are shown in 
figures~\ref{mm1}-\ref{mm5} in the back of the paper. 
Similar to modifications (a)-(e) of generating morphisms, there  
are modifications (a)-(e) of movie moves and they should be included 
in the list. See \cite{BaezLang} for 
a detailed discussion.  

Move 31 is not shown. It says that given horizontally composable 
2-morphisms $\alpha: f\Rightarrow f'$ and $\beta: g\Rightarrow g',$ 
there is an equality  
$(\alpha\cdot \Id)(\Id \cdot \beta) =(\Id \cdot \beta)(\alpha\cdot \Id)$
of 2-morphisms from $fg$ to $f'g'.$   

Figures~\ref{mm1}-\ref{mm3} show local moves, while figures~\ref{mm4}, 
\ref{mm5} show semi-local moves. Little squares in semi-local moves could be 
U-turns or crossings.

\section{Bimodule homomorphisms} 

For a ring $A$ denote by $\mc{K}(A)$ the category of bounded complexes 
of $A$-bimodules up to homotopies of complexes.The objects of $\mc{K}(A)$ 
are bounded complexes of $A$-bimodules, the morphisms are morphisms of 
complexes of bimodules quotiented by homotopic to $0$ morphisms. 

We say 
that a complex of bimodules $M \in \mc{K}(A)$ is \emph{invertible} if there 
exists $N \in \mc{K}(A)$ such that $N\otimes_A M\cong A$ and 
$M\otimes_A N\cong A$ 
in $\mc{K}(A).$ Here $A$ denotes 
the complex $0\lra A\lra 0$ with $A$ in cohomological degree $0$ and 
the usual left and right multiplication action of $A$ on itself. 

Let $Z(A)$ be the center of $A.$ 

\begin{prop} If $M$ is invertible then 
\begin{equation*}
  \mathrm{Hom}_{\mc{K}(A)}(M,M)\cong
  \mathrm{Hom}_{A\otimes A^o}(A,A)\cong Z(A), 
\end{equation*}   
\end{prop} 

\emph{Proof:} 
The second isomorphism is obvious, since endomorphisms of $A$ as an 
$A$-bimodule are multiplications by central elements. 

Consider the following sequence of ring homomorphisms: 
\begin{eqnarray*}
 & &  \mathrm{End}_{\mc{K}(A)}(M) \stackrel{f}{\lra}
  \mathrm{End}_{\mc{K}(A)}(M\otimes_A N) \stackrel{g}{\lra} 
   \mathrm{End}_{\mc{K}(A)}((M\otimes_A N)\otimes_A M ) \cong    \\
 & &   \mathrm{End}_{\mc{K}(A)}(M\otimes_A (N\otimes_A M) ) \cong  
      \mathrm{End}_{\mc{K}(A)}(M \otimes_A   A  ) \cong 
  \mathrm{End}_{\mc{K}(A)}(M), 
  \end{eqnarray*}
where $f,$ respectively $g,$ is tensoring with the 
identity endomorphism of $N,$ respectively $M.$ The composition 
$gf$ is the identity, thus $f$ is injective.
Multiplication on the left by central elements makes 
each of the above rings a $Z(A)$-module, and $f$ and $g$ are $Z(A)$-module 
homomorphisms. $f$ and $g$ take identity endomorphisms to identity 
endomorphisms, and  $\mathrm{End}_{\mc{K}(A)}(M\otimes_A N)\cong 
\mathrm{End}_{\mc{K}(A)}(A) = Z(A).$ Therefore, $f$ is surjective, since 
  $f(\mathrm{id})= \mathrm{id}$ generates 
  $\mathrm{End}_{\mc{K}(A)}(M\otimes_A N)$
  as a $Z(A)$-module. Thus, $f$ and $g$ are isomorphisms. $\square$

If $A$ is graded, denote by $\mc{K}(A)$ the category of 
bounded complexes of graded $A$-bimodules (with grading-preserving 
differential) up to homotopies. The morphisms 
are grading-preserving homomorphisms of complexes (modulo homotopies). 
If $M$ is an invertible complex in $\mc{K}(A)$ then 
$\mathrm{Hom}_{\mc{K}(A)}(M,M) \cong Z_0(A),$ the degree $0$ component 
of the center of $A.$ 
Furthermore, the group of automorphisms of $M$ in $\mc{K}(A)$ is 
isomorphic to $Z_0^{\ast}(A),$ the group of invertible elements in $Z_0(A).$ 

\vsp 

We now specialize to the rings $H^n.$ 

\begin{prop} The only invertible central elements of degree $0$ in $H^n$ 
are $\pm 1:$ 
\begin{equation*}
Z_0^{\ast}(H^n) \cong \{ \pm 1 \}.
\end{equation*} 
\end{prop} 

\emph{Proof:} A degree $0$ element 
of $H^n$ has the form $v=\sum v_a e_a$ where $e_a$ is the minimal 
idempotent corresponding to the crossingless matching $a$ and $v_a\in \Z.$ 
For any $a,b$ choose $x\in \hsm _a(H^n)_b, x\not= 0.$ Then 
$vx = v_a x$ and $xv = v_b x.$ Therefore, if $v$ is central, $v_a=v_b$ for 
all $a,b,$ and  $v=m \sum e_a = m,$ for some integer $m,$ so that 
$Z_0(H^n) \cong \Z.$ The proposition follows. $\square$ 

\emph{Remark:} We investigated the center of $H^n$ 
(and not just its degree $0$ component) in \cite{me:springer}. 
It turned out to be 
isomorphic to the cohomology ring of the $(n,n)$ Springer fiber. 

\begin{corollary} \label{inv-cor}
If $M$ is an invertible complex of graded $H^n$-bimodules, 
then $\mathrm{Id}$ and $-\mathrm{Id}$ are the only degree $0$ 
automorphisms of $M.$ 
\end{corollary}

We use notation $\cK_n^m$ from \cite{me:tangles} for the category 
of bounded complexes of geometric $(H^m,H^n)$-bimodules up to chain 
homotopies. A bimodule is geometric if it is isomorphic to a finite 
direct sum of bimodules $\cF(a)\{ i\},$ for flat tangles $a$ and 
$i\in \Z$ (recall that $\{ i \}$ denotes shift in the grading by $i$).  
Morphisms in $\cK_n^m$ are grading-preserving homomorphisms 
of complexes up to chain homotopies.

From Corollary~\ref{inv-cor} we derive 

\begin{corollary} If $f: M\to N$ is an isomorphism of invertible 
objects in $\cK_n^n$ then the only other isomorphism of $M$ and $N$ is 
$-f.$ \label{only-minus}
\end{corollary}

\vsp 

For now on we assume that the reader is familiar with
the construction of \cite[Section 2]{me:tangles}, which to a surface $S$
embedded in $\R^3$ and viewed as a cobordism between flat tangles $a$ and 
$b$ assigns a bimodule homomorphism $\cF(S): \cF(a)\lra \cF(b).$ 

Standard cobordisms in $\R^3$ (the first is a birth move, the second and 
third are saddle points, the fourth is a death move)

\begin{eqnarray*}
 \mathrm{Vert}_{n-1} & \Longrightarrow & \cap_{i,n}\cup_{i,n-1}, \\
 \cup_{i,n-1}\cap_{i,n} & \Longrightarrow & \mathrm{Vert}_n,   \\
 \mathrm{Vert}_n & \Longrightarrow & \cup_{i,n-1} \cap_{i,n}, \\ 
 \cap_{i,n}\cup_{i,n-1} & \Longrightarrow & \mathrm{Vert}_{n-1}
\end{eqnarray*}
induce grading-preserving bimodule homomorphisms
\begin{eqnarray} 
 \label{eqn-1}
   & & H^{n-1} \lra   \cF(\cap_{i,n})\otimes_{H^n} \cF(\cup_{i,n-1})\{ 1\}, \\
 \label{eqn-2}
  & & \cF(\cup_{i,n-1})\otimes_{H^{n-1}} \cF(\cap_{i,n})\{ 1\}\lra H^n, \\
 \label{eqn-3}
  & & H^{n} \lra\cF(\cup_{i,n-1})\otimes_{H^{n-1}}\cF(\cap_{i,n})\{ -1\}, \\ 
 \label{eqn-4}
   & & \cF(\cap_{i,n})\otimes_{H^{n}} \cF(\cup_{i,n-1})\{ -1\} \lra H^{n-1},
\end{eqnarray}
(we used that $\cF(\mathrm{Vert}_n)\cong H^n, 
\cF(\cup_{i,n-1} \cap_{i,n}) \cong 
\cF(\cup_{i,n-1})\otimes_{H^{n-1}}\cF(\cap_{i,n}),$ etc.) 
Isotopies between compositions of these cobordisms translate into 
relations between homomorphisms. These relations imply that 
the functors  of tensoring 
with $\cF(\cap_{i,n})$ and $\cF(\cup_{i,n-1})$ are biadjoint, up to 
grading shifts. Precisely, 
let $F_{\cup}$ be the functor of tensoring with $\cF(\cup_{i,n-1})$ and 
$F_{\cap}$ the functor of tensoring with $\cF(\cap_{i,n})$ (viewed 
as functors between categories of $H^n$ and $H^{n-1}$-modules).  

\begin{prop} \label{biadj} 
$F_{\cup}\{ 1\}$ is left adjoint to $F_{\cap},$ and 
  $F_{\cap}\{ -1\}$ is left adjoint to $F_{\cup}.$ 
\end{prop} 

\begin{corollary} The only grading-preserving endomorphisms of 
bimodules $\cF(\cap_{i,n})$ and $\cF(\cup_{i,n-1})$ are multiplications 
by integers. The only grading-preserving automorphisms of bimodules 
$\cF(\cap_{i,n})$ and 
$\cF(\cup_{i,n-1})$ are $\Id$ and $-\Id.$ Moreover, these bimodules 
have no graded endomorphisms of negative degree. 
\end{corollary}

\emph{Proof  of corollary:} From adjointness, 
\begin{eqnarray*}
 \Hom_{(n,n-1)}(\cF(\cup_{i,n-1}), \cF(\cup_{i,n-1})) & \cong &  
  \Hom_{(n-1,n-1)} ( H^{n-1}, \cF(\cap_{i,n}\cup_{i,n-1}) \{ 1\} ) \\
   & \cong &  
  \Hom_{(n-1,n-1)} ( H^{n-1}, H^{n-1}\oplus H^{n-1}\{ 2\})    \\
  & \cong & \Hom_{(n-1,n-1)} ( H^{n-1}, H^{n-1})   \\
  & \cong & \Z. \\ 
\end{eqnarray*}
Subscripts of the form $(m,n)$ in the above formula mean that  
the homomorphisms considered are those of graded $(H^m,H^n)$-bimodules. 
We used that 
\begin{equation*}
 \Hom_{(n-1,n-1)} ( H^{n-1}, H^{n-1}\{ k\})=0
 \end{equation*}
for any positive $k,$ since the ring $H^{n-1}$ is $\Z_+$-graded. 
Similar computations establish the result for $\cF(\cap_{i,n})$ and the last 
claim of the corollary. $\square$ 

\begin{corollary} \label{tens-cap}
If $M$ is a tensor product of $\cF(\cap_{i,n})$ and 
 invertible complexes of bimodules, and if $f: M\to N$ is an isomorphism in 
$\cK_n^{n-1},$ then the only other isomorphism from $M$ to $N$ is $-f.$ 
Same with $\cup_{i,n-1}$ instead of $\cap_{i,n}.$ 
\end{corollary}

More generally, let $b$ be a flat tangle without closed components (circles), 
with $k$ arcs connecting 
bottom endpoints, $l$ arcs connecting top endpoints, and some 
number of arcs connecting a top endpoint with a bottom endpoint. Let $W(b)$ be 
the reflection of $b$ about the $x$-axis. Representing 
$b$ as a product of $U$-turns and using Proposition~\ref{biadj} 
repeatedly we obtain 

\begin{prop} The functor of tensoring with the bimodule 
$\cF(W(b)) \{ k-l\}$ is left adjoint to tensoring with $\cF(b)$ and the 
functor of tensoring 
with $\cF(W(b)) \{ l-k\}$ is right adjoint to tensoring with $\cF(b).$ 
\end{prop} 

\begin{corollary} If $b$ is a flat tangle without closed components then 
the only grading-preserving endomorphisms of the bimodule $\cF(b)$ are 
multiplications by integers, the only grading-preserving automorphisms 
are $\Id$ and $-\Id,$ and $\cF(b)$ has no endomorphisms of negative 
degree. If $f: M \to \cF(b)$ is an isomorphism in $\cK_n^m$ then the 
only other isomorphism between $M$ and $\cF(b)$ is $-f.$ \label{flat-minus}
\end{corollary}

\section{The 2-functor}

We introduce two 2-categories $\mathbb{K}$ and $\widehat{\mathbb{K}}.$ 

Objects of $\mathbb{K}$ are non-negative integers, 1-morphisms from $n$ 
to $m$ are objects of $\cK_n^m,$ and 2-morphisms between $M,N\in \cK_n^m$ 
are $\Hom_{\cK_n^m}(M,N),$ grading-preserving 
morphisms of complexes of bimodules up to 
chain homotopies. Composition of 1-morphisms is given by tensor product 
of complexes.  

The 2-category $\widehat{\mathbb{K}}$ has the same objects and 1-morphisms 
as $\mathbb{K}$ but the 2-morphisms are 
\begin{equation*}
  \Hom_{\widehat{\mathbb{K}}}(M,N)\define \oplusop{i\in \Z}
      \Hom_{\cK_n^m}(M,N \{i\}) /\{ \pm 1\}, 
\end{equation*}
that is, the morphisms are all homomorphisms (not just grading-preserving), 
and each homomorphism is identified with its negative. The set of 
2-morphisms between two 1-morphisms is no longer an abelian group.  

\vsp 

We next construct a 2-functor $\cF: \cC \to \widehat{\mathbb{K}}.$ 
This 2-functor takes object $s$ of $\cC$ to the 
object $|s|$ of $\widehat{\mathbb{K}}.$ It takes generating 1-morphisms 
of $\cC$ to complexes of bimodules in the same way as in 
\cite[Sections 2.7, 3.4]{me:tangles}. Recall that a U-turn $b$ (and, more 
generally, any flat tangle) is taken to the complex 
$0 \lra \cF(b) \lra 0$ where $\cF(b)$ is the bimodule associated to $b.$ 
A crossing $r$ gives rise to its two resolutions $r(0)$ and $r(1)$ and 
a grading-preserving bimodule map $\psi: \cF(r(0))\lra \cF(r(1))\{-1\}.$  
Then $\cF(r)$ is defined as the complex 
\begin{equation*} 
 0 \lra \cF(r(0))\stackrel{\psi}{\lra} \cF(r(1))\{-1\} \lra 0
\end{equation*}
with a suitable grading shift computed from the orientation of $r$ 
near its crossing. 

The 2-functor $\cF$ takes composition of 1-morphisms 
to the tensor product of complexes: 

\begin{equation*}
\cF(ab) \define \cF(a) \otimes_{H^n} \cF(b), 
\end{equation*}
where $a,$ respectively $b,$ has $2n$ bottom, respectively, $2n$ top 
endpoints. 

To a Reidemeister move between diagrams $a$ and $b$ (see 
figure~\ref{gen1mor}) we assign an isomorphism of bimodule complexes 
$\cF(a) \stackrel{\cong}{\lra} \cF(b)$ constructed 
in \cite[Section 4]{me:tangles}. The Reidemeister III move 
has several versions, depending on the directions of overcrossings. 
In \cite{me:tangles} we described an isomorphism in $\cK_n^n$ between 
complexes $\cF(a)$ and $\cF(b)$ for only one version of this move. 
Other versions can be expressed via compositions of this version 
with isotopies and type II moves. The compositions induce isomorphisms 
between $\cF(a)$ and $\cF(b)$ which we assign to these other version 
of the Reidemeister III move.  Note that, since Reidemeister move 
diagrams in figure~\ref{gen1mor} are either braids or 
composition of a U-turn and a braid, a grading-preserving isomorphism 
between $\cF(a)$ and $\cF(b)$ 
is unique up to minus sign, by Corollaries~\ref{only-minus}, \ref{tens-cap}.

The diagrams of birth, death, saddle point, and T-move do not involve 
crossings. These movies can be viewed as presentations of surfaces 
embedded in $\R^3,$ and to them we assign bimodule homomorphisms using the 
construction of Proposition 5 of \cite{me:tangles}.  
To the birth 2-morphism we assign the unit map $\iota: \Z \to \cA,$ to the 
death 2-morphism the counit map $\epsilon: \cA \to \Z.$ More precisely, 
birth and death moves happen inside $\Vertical_{n}$ diagrams, and 
the maps are 
\begin{equation*}
\iota\otimes \Id_{H^n}: H^n \lra \cA\otimes_{\Z} H^n, \hspace{0.15in}
 \epsilon\otimes \Id_{H^n}: \cA\otimes_{\Z} H^n \lra H^n. 
\end{equation*}

To the diagrams of saddle point cobordisms between $\Vertical_{n}$ and 
$\cup_{i,n-1}\cap_{i,n}$ we assign bimodule maps (\ref{eqn-2}), (\ref{eqn-3}). 
These maps are, up to sign, 
the only maps of degree $1$ between $\cF(\Vertical_n)\cong H^n$ 
and $\cF(\cup_{i,n-1}\cap_{i,n})$ 
that generate the abelian group (isomorphic to $\Z$) of all degree $1$ 
homomorphisms between these bimodules. 

The natural isotopy between the two diagrams $a,b$ in the H-move  
(figure~\ref{2mor2}) induces 
an isomorphisms of complexes $\cF(a)\cong \cF(b),$ which we assign 
to this 2-morphism.  

A frame of an N-move between diagrams $a$ and $b$ 
contains two little squares, and each square is either a U-turn or a crossing. 
The N-move is an isotopy from $a$ to $b.$ This isotopy induces a 
canonical isomorphism of complexes $\cF(a) \stackrel{\cong}{\lra}\cF(b)$
(see \cite[Sections 4.1, 4.2]{me:tangles}), which we assign to the 
N-move.   

\begin{theorem} The above correspondence extends uniquely to a 2-functor 
 \begin{equation*}
 \cF: \mc{C}\to \widehat{\mathbb{K}}.
 \end{equation*}
\end{theorem}

This theorem is proved in Section~\ref{proof}. 

$\mc{C}$ is a combinatorial 
realization of the 2-category $\mc{T}$ of tangle cobordisms, that is, 
the natural 2-functor $\mc{C}\to \mc{T}$ is an equivalence of 2-categories, 
see \cite{BaezLang}. This result is also valid for oriented tangles. 

As a corollary, we obtain a 2-functor, also denoted $\cF,$ from the 
2-category $\mc{T}$ of even unframed oriented tangle cobordisms to
$\widehat{\mathbb{K}}.$ The homomorphism of complexes of graded 
 $(H^m,H^n)$-bimodules 
assigned to the cobordism $S$ between $(m,n)$-tangles has degree 
$n+m-\chi(S),$ where $\chi(S)$ is the Euler characteristic of $S.$

\section{Proof}  \label{proof} 

When looking at a particular movie move 
we denote the top frame by $b_1,$ the bottom frame by $b_2,$ the left movie 
by $S_l$ and the right movie by $S_r.$ 

Movies $S_l$ and $S_r$ induce homomorphisms $\cF(S_l)$ and 
$\cF(S_r)$ from $\cF(b_1)$ to $\cF(b_2).$ We need to show that 
$\cF(S_l) = \pm \cF(S_r)$ in $\cK_n^m.$ 

\vspace{0.1in}

Moves 1, 2, 3, 4, 5 say that composing a Reidemeister move with 
its inverse is equivalent to doing nothing. The isomorphism in 
$\cK_n^m$ assigned to the inverse of a Reidemeister move equals the 
inverse of the isomorphism assigned to the move. Therefore, 
 $\cF(S_l)= \cF(S_r)$ for each of these moves.  

\vspace{0.1in}

Movies $S_l$ and $S_r$ in move 6 consist of Reidemeister moves and 
relative height shifts of distant crossings. 
 The complexes $\cF(b_1)$ and $\cF(b_2)$ are invertible (since $b_1$ and 
$b_2$ are braids), and 
\begin{equation*}
\cF(S_l), \cF(S_r): \cF(b_1)\lra \cF(b_2)
\end{equation*}
 are two isomorphisms of these complexes in $\cK_n^n.$ By 
Corollary~\ref{only-minus}, either  $\cF(S_l) = \cF(S_r)$ or 
 $\cF(S_l) = -\cF(S_r).$

This proof works simultaneously for all versions of move 6. Identical 
 argument takes care of moves 12, 13,  23a and  25 (and of moves 
 3, 4, 5  as well, although the latter have already been dealt with).  
 
\vspace{0.1in} 

Each movie in move 7 is a composition of Reidemeister moves, 
thus, $\cF(S_l)$ and $\cF(S_r)$ are grading-preserving isomorphisms 
(in $\cK_n^{n-1}$) of 
complexes $\cF(b_1)$ and $\cF(b_2).$ Since $b_1$ and $b_2$ are 
given by composing $\cap_{i,n}$ with braids, $\cF(b_1)$ and $\cF(b_2)$ 
are tensor products of $\cF(\cap_{i,n})$ with invertible complexes (the 
index $i$ is different for $b_1$ and $b_2$). 
By Corollary~\ref{tens-cap}, $\cF(S_l)$ differs from  $\cF(S_r)$ by 
at most a minus sign. Other versions of this move follow suit. 

Identical arguments takes care of moves 11, 14, and 26. 

\vspace{0.1in}

Moves 8, 9, 10, 23b, 24 do not involve any crossings and the invariance 
of $\cF$ follows from 
Proposition 6 of \cite{me:tangles}, since these moves are saying that 
certain surfaces in $\R^3$ are isotopic.  

\vspace{0.1in}

Both movies in move 21 consist of isotopies and a Reidemeister move. 
Therefore, $\cF(S_l)$ and $\cF(S_r)$ are isomorphisms in $\cK_n^n.$ 
The bottom diagram $b_2$ is a flat tangle without
 closed components (circles).  
Corollary~\ref{flat-minus} implies that any two isomorphisms from 
$\cF(b_1)$ to $\cF(b_2)$ differ by sign at most. Similar arguments 
take care of moves 15-20 (use Corollary~\ref{flat-minus} and 
its generalization from $b$ to $b_1 b b_2$ where $b_1$ and $b_2$ 
are braids). Alternatively, the invariance of $\cF$ under semi-local 
moves 15-20, 22 follows by observing that height shifts of U-turns 
and crossings don't do anything to our complexes of bimodules and 
maps between them. 

\vspace{0.1in}

The first frame change in both movies in move 28 is birth, which is 
then followed by a Reidemeister move and an isotopy (H-move). 
Decompose $S_l = R_l Q_l$ and $S_r = R_r Q_r$ where $Q_l, Q_r$ are births. 
Denote by $b'_l, b'_r$ second frames from the top in the left 
and right movies. Note that $\cF(b'_l)\cong \cF(b'_r)\cong \cA\otimes H^n,$ 
and $\cF(R_r)^{-1}\cF(R_l): \cF(b'_l) \to \cF(b'_r)$ 
is a grading-preserving isomorphism, while 
$\cF(Q_l), \cF(Q_r): H^n \lra   \cA\otimes_{\Z} H^n$
have degree $-1.$ Both $\cF(Q_r)$ and $ \cF(R_r)^{-1}\cF(S_l)$ generate 
the abelian group $\Z$ of degree $-1$ homomorphisms from $H^n\cong \cF(b_1)$ 
to $\cA\otimes_{\Z} H^n \cong \cF(b'_r).$ Therefore,  
$\cF(Q_r)$ and $ \cF(R_r)^{-1}\cF(S_l)$ differ by 
at most minus sign, and $\cF(S_l), \cF(S_r)$ differ by at most minus 
sign. 

Invariance of $\pm \cF$ under moves 22 and 27 follows from similar arguments. 

\vspace{0.1in}

Both movies in move 29 consist of a Reidemeister move followed by a
saddle point 2-morphism. Both movies induce degree 1 homomorphisms 
from $\cF(b_1)$ to $\cF(b_2).$ Since the homomorphism assigned to 
the saddle point generates the abelian group (isomorphic to $\Z$) 
of degree 1 homomorphisms 
from $\cF(\cup_{i,n-1}\cap_{i,n})$ to $H^n\cong \cF(\Vertical_n),$
 we see that both $\cF(S_l)$ and $\cF(S_r)$ are generators of 
\begin{equation*}
\Hom_{\cK_n^n}(\cF(b_1)\{ 1\}, \cF(b_2))\cong \Z,
\end{equation*}    
and differ by  at most minus sign. Move 30 follows similarly.  

\vsp 

Given a ring $A$ and homomorphisms $f_1: M_1\to N_1,$ resp.  
$f_2: M_2\to N_2$ of complexes of right, resp. left, $A$-modules, the map 
  \begin{equation*}
   f_1\otimes f_2: M_1\otimes_A M_2\lra N_1\otimes_A N_2
  \end{equation*}
 can be written in two ways: 
as $(f_1\otimes \Id)(\Id \otimes f_2)$ and as 
$(\Id\otimes f_2)(f_1\otimes \Id).$ This observation takes care of 
move 31.   $\square$

\newpage

 \newpage 
 
 \newpage 
 
 \newpage 
 
\newpage 
 
\end{document}